\catcode`\ =11%
\def\space{ }%
\catcode`\ =10%
\def\begincolor#1{}%

\def\endcolor{}%

\def\colored#1#2{\begincolor{#1}#2\endcolor}%

\def\linkcolor{[0.0 0.0 1.0]}%

\def\setlink#1{\colored{\linkcolor}{#1}}%

\def\link#1#2{\setbox0\hbox{\setlink{#1}}%
   \leavevmode\box0\relax}%
\def\dest#1{\hbox{\raise14pt\hbox{}}}%

\def\bookm#1#2{}
\input amssym.tex
\font\twbf=cmbx12

\font\fiverm=cmr5
\font\sc=cmcsc10

\def\bg{\bigskip \goodbreak}
\def\d{{\rm d}}
\def\ds{\displaystyle}
\def\abs#1{\left\vert\, #1\,\right\vert }
\def\floor#1{\left\lfloor #1\right\rfloor }
\def\ceiling#1{\left\lceil #1\right\rceil }
\def\rond{{\scriptstyle\circ}}
\def\ci{{\rm i}}
\def\cqfd{\hfill\hbox{\vrule\vbox to 4pt{\hrule width 4pt\vfill \hrule}\vrule}}
\def\nat{\Bbb{N}}
\def\real{\Bbb{R}}
\def\frac#1#2{{#1\over#2}}
\def\ref#1&#2&#3&#4&#5{
\item{[{\bf\ignorespaces#1\unskip}]}
{\sc\ignorespaces#2\unskip},\
{\rm\ignorespaces#3\unskip}\
{\sl\ignorespaces#4\unskip\/}\
{\rm\ignorespaces#5\unskip}}

%
%

\centerline{\twbf  Inequalities related to the error function}\par
\bg
\centerline{{\sc Omran Kouba}}\par
\centerline{Department of Mathematics}\par
\centerline{{\sl Higher Institute for Applied Sciences and Technology}}\par
\centerline{P.O. Box 31983, Damascus, Syria.}\par
\centerline{{\it E-mail} : omran\_kouba@hiast.edu.sy}\par
\bg
\par
\parindent=50pt
{\narrower\smallskip\noindent {\bf Abstract:} 
In this note we consider inequalities involving the function
$$ \varphi:\Bbb{R}\rightarrow\Bbb{R},~\varphi(x)=e^{x^2/2}\int_x^\infty e^{-t^2/2}\d t.$$
Our methodes give new proofs of some known inequalities of Komatsu [\link{1}{r1}], 
and Szarek and Werner [\link{2}{r2}], and also produce two families of inequalities 
that give upper and lower bounds for $\varphi$. Moreover, the continued fractions expansion
of $\varphi(x)$ for $x>0$ is obtained.\smallskip}
\parindent=0pt
\bg
\bookm{1}{Introduction}
\noindent{\bf 1. Introduction.}
\par
Komatsu's Inequalities [\link{1}{r1}] assert that 
$$\forall x\in\Bbb{R}_+,\qquad{2\over x+\sqrt{x^2+4}}<\varphi(x)<{2\over x+\sqrt{x^2+2}}$$
In fact, the inequality on the left is true for every real $x$.
In  [\link{2}{r2}] the authors sharpened the inequality on the right to become
$$\forall x>-1,\qquad \varphi(x)<{4\over 3x+\sqrt{x^2+8}}.$$
\par
In this note we present two families of inequalities that  give upper and lower bounds for $\varphi$, and give the above-mentioned inequalities as special cases.
\bg
\bookm{1}{Preliminaries}
\noindent{\bf 2. Preliminaries.}
\par
In the following lemma, we introduce two sequences of polynomials that will play an
important role in what follows.\par
\bookm{2}{Definition of Pn and Qn}\dest{p1}
\proclaim Lemma 1. For every non negative integer $n$, there exists a unique couple
 $(P_n,Q_n)$ of polynomials that satisfy
$$\forall\, x\in\real,\quad\varphi^{(n)}(x)=P_n(x)\varphi(x)-Q_n(x)\eqno{(1)}\dest{e1}$$
Moreover, these polynomials are defined, starting from $(P_0,Q_0)=(1,0)$, by the recurrence relations
$$\forall\, n\in\nat,\quad P_{n+1}=XP_n+P^\prime_n\eqno{(2)}\dest{e2}$$
$$\forall\, n\in\nat,\quad Q_{n+1}= P_n+Q^\prime_n\eqno{(3)}\dest{e3}$$
\par
\noindent{\bf Proof.}\par
The statement concerning uniqueness is obvious, since a property of the form
 $$\forall\, x\in\real,\quad U(x)\varphi(x)+V(x)=0$$ with  two polynomials $U$ and $V$ would imply $\lim_{x\to-\infty}U(x)=0$ 
because $\varphi(x)$ is equivalent to $\sqrt{2\pi}e^{x^2/2}$ in a neighbourhood of $-\infty$.\par
On the other hand, it is clear that $\varphi^\prime(x)=x\varphi(x)-1$ which proves\link{ $(1)$}{e1}
for $n=1$ with $(P_1,Q_1)=(X,1)$. Assuming the existence of two polynomials $P_n$ and $Q_n$ satisfying
 \link{$(1)$}{e1} we conclude by differentiation that
$$\eqalign{\varphi^{(n+1)}=&(P_n\varphi-Q_n)^\prime=P_n\varphi^\prime+P_n^\prime\varphi-Q^\prime_n\cr
=&(X\varphi-1)P_n+P_n^\prime\varphi-Q^\prime_n\cr
=&(\underbrace{XP_n+P_n^\prime}_{P_{n+1}})\varphi-(\underbrace{P_n+Q^\prime_n}_{Q_{n+1}})}$$
and this achieves the proof by mathematical induction.\cqfd 
\par
\bg
\bookm{1}{Properties of the polynomials Pn and Qn}
\noindent{\bf 3. Properties of the sequences $(P_n)_{n\in\nat}$ and $(Q_n)_{n\in\nat}$.}
\bg
It is clear, using the recurrence relations \link{$(2)$}{e2} and \link{$(3)$}{e3} that the coefficients
 of $P_n$ and $Q_n$ are non-negative integers and that the dominant monomial in $P_n$ is $X^n$, and the dominant monomial in 
$Q_n$ is $X^{n-1}$, for $n\geq1$.
In the following proposition we summarize some elementary properties of these polynomials.\bg
\bookm{2}{Recurrence formulae}
\parindent=20pt
\proclaim Proposition 2.  The sequences $(P_n)_{n\in\nat}$ and $(Q_n)_{n\in\nat}$ satisfy the following properties.
\item{$1^\rond.$} $(P_0,P_1)=(1,X)$, and for all $n\geq1$ we have $P_{n+1}=XP_n+nP_{n-1}$.\dest{e4}\hfill$(4)$
\item{$2^\rond.$} $(Q_0,Q_1)=(0,1)$, and for all $n\geq1$ we have $Q_{n+1}=XQ_n+nQ_{n-1}$.\dest{e5}\hfill$(5)$
\item{$3^\rond.$} For all $n\geq1$ we have $P^\prime_{n}=nP_{n-1}$.\dest{e6}\hfill$(6)$
\par
\parindent=0pt
{\bf Proof.}\par
Using the fact that $\varphi^\prime(x)=x\varphi(x)-1$, we can write, for every real $x$ and every $n\geq 1$,
$$\varphi^{(n+1)}(x)=(x\varphi(x)-1)^{(n)}=x\varphi^{(n)}(x)+n\varphi^{(n-1)}(x)$$
that is
$$P_{n+1}\varphi-Q_{n+1}=(XP_n+nP_{n-1})\varphi-(XQ_n+nQ_{n-1})$$
and the uniqueness statement in \link{lemma 1}{p1} implies \link{$(4)$}{e4} and \link{$(5)$}{e5}. \par
Finally, comparing  \link{$(4)$}{e4}  and  \link{$(2)$}{e2} proves \link{$(6)$}{e6}.\cqfd
\bg
The next proposition contains two representations of the polynomials $(P_n)_{n\in\nat}$.\bg
\bookm{2}{Integral representation of P_n}
\parindent=20pt
\proclaim Proposition 3.  The sequence $(P_n)_{n\in\nat}$ satisfies the following properties.
\item{$1^\rond.$}  For all $n\geq0$, and every real $x$,  we have $P_{n}(x)=e^{-x^2/2}\left(e^{x^2/2}\right)^{(n)}$.\dest{e7}\hfill$(7)$
\item{$2^\rond.$}  For all $n\geq0$, and every real $x$, we have $$P_{n}(x)={1\over\sqrt{2\pi}}\int_\real t^n\exp\left(-{(t-x)^2\over2}\right)\d t\dest{e8}\eqno{(8)}$$
\par
\parindent=0pt
{\bf Proof.}\par
{$1^\rond.$} Let us define $$h_n:\real\rightarrow\real, h_n(x)=e^{-x^2/2}\left(e^{x^2/2}\right)^{(n)}.$$ Clearly, we have
$h_0\equiv 1\equiv P_0$ and $\forall\,x\in\real,~h_1(x)=x=P_1(x)$. 
Moreover, for all $n\geq1$, and every real $x$, we have
$$\left(e^{x^2/2}\right)^{(n+1)}=\left(xe^{x^2/2}\right)^{(n)}=x\left(e^{x^2/2}\right)^{(n)}+n\left(e^{x^2/2}\right)^{(n-1)}$$
that is
$$\forall\,n\in\nat^*,~\forall\,x\in\real,\quad h_{n+1}(x)=xh_{n}(x)+nh_{n-1}(x)$$
So, the sequence $(h_n)_{n\in\nat}$ satisfies the recurrence relation \link{$(4)$}{e4} with the same initial conditions as the sequence
 $(P_n)_{n\in\nat}$. This proves that  $\forall\,n\in\nat,~P_n=h_n$, that is  \link{$(7)$}{e7}.\par
{$2^\rond.$} Using the change of variable $u\leftarrow t-x$ in the integral ${1\over\sqrt{2\pi}}\int_\real e^{-u^2/2}\d u=1$ we find that
$$\forall\,x\in\real,\quad e^{x^2/2}={1\over\sqrt{2\pi}}\int_\real e^{t x} e^{-t^2/2}\d t.$$
And since differentiation under integral sign is legitimate, we conclude, using \link{$(7)$}{e7}, that  
$$\forall\,x\in\real,\quad P_n(x)={1\over\sqrt{2\pi}}\int_\real t^ne^{-(t- x)^2/2}\d t={1\over\sqrt{2\pi}}\int_\real( u+x)^ne^{-u^2/2}\d u.\eqno\cqfd$$
\bg
\bookm{2}{Parity and lower bounds}
\proclaim Corollary 4.  The polynomials $(P_n)_{n\in\nat}$ are given by the following formula.
 $$P_{n}=\sum_{k=0}^{\floor{n/2}}{n!\over 2^kk!(n-2k)!}X^{n-2k}\dest{e9}\eqno{(9)}$$
In particular, $P_n$ has the same parity as $n$, and
$$\forall\,x\in\real,\quad P_{2n}(x)\geq{(2n)!\over 2^n n!}\quad\hbox{and}\quad
{P_{2n+1}(x)\over x}\geq{(2n+1)!\over 2^n n!}\dest{e10}\eqno{(10)}$$
\par
{\bf Proof.}\par
In fact, for every non-negative integer $n$, and every real $x$, we have
$$\eqalign{P_n(x)=&{1\over\sqrt{2\pi}}\int_\real( u+x)^ne^{-u^2/2}\d u\cr
=&\sum_{k=0}^n\left({n\atop k}\right)x^{n-k}{1\over\sqrt{2\pi}}\int_\real u^ke^{-u^2/2}\d u\cr
=&\sum_{0\leq2k\leq n}^n\left({n\atop 2k}\right)x^{n-2k}{1\over\sqrt{2\pi}}\int_\real u^{2k}e^{-u^2/2}\d u}$$
but we know that 
$${1\over\sqrt{2\pi}}\int_\real u^{2k}e^{-u^2/2}\d u=2^k{\Gamma(k+{1\over2})\over\sqrt\pi}={(2k)!\over 2^kk!}$$
which proves
$$P_n(x)=\sum_{0\leq2k\leq n}\left({n\atop 2k}\right){(2k)!\over 2^kk!}x^{n-2k}$$
and this is the desired result.\cqfd
\par
\bg
\bookm{2}{Link with Hermite polynomials}
{\bf Remark.} It is easy to see that $P_n$ is linked to the well known {\sc Hermite} polynomials usually denoted as $H_n$ in the following way:
$$P_n(X)=\left({-\ci\over\sqrt2}\right)H_n\left({\ci\over\sqrt2}X\right)$$
\bg
Expressions for $Q_n$ are not so simple, as we will see in the following proposition.\bg
\bookm{2}{Qn in terms of Pn}
\parindent=20pt
\proclaim Proposition 5.  The sequence $(Q_n)_{n\in\nat}$ satisfies the following properties.
\item{$1^\rond.$}  For all $n\geq0$, we have 
$$Q_{n+1}=\sum_{k=0}^{\floor{n/2}}{(n-k)!\over(n-2k)!}P_{n-2k}\dest{e11}\eqno{(11)}$$
\item{$2^\rond.$}  For all $n\geq0$, we have 
$$Q_{n+1}=\sum_{k=0}^{\floor{n/2}}{1\over(n-2k)!}\left(\sum_{j=0}^k{(n-k+j)!\over2^jj!}\right)X^{n-2k}\dest{e12}\eqno{(12)}$$
\par
\parindent=0pt
{\bf Proof.}\par
{$1^\rond.$} Consider $(R_n)_{n\geq 1}$ defined by 
$$\forall\,n\in\nat^*,\quad R_{n+1}=\sum_{k\geq0}P_{n-k}^{(k)}.$$
It is straightforward to see that for $n\geq1$ we have
$$\eqalign{R_{n+1}=&P_n+\sum_{k\geq1}P_{n-k}^{(k)}=P_n+\sum_{k\geq0}P_{n-1-k}^{(k+1)}\cr
=&P_n+\left(\sum_{k\geq0}P_{n-1-k}^{(k)}\right)^\prime=P_n+R_n^\prime}$$
Using this, the fact that $R_1=Q_1=1$, and formula  \link{$(3)$}{e3}, 
we conclude by induction that $R_n=Q_n$ for every $n\geq1$. Moreover, 
using \link{$(6)$}{e6}, we can write :
$$P_{n-k}^{(k)}=(n-k)(n-1-k)\cdots(n-2k+1)P_{n-2k}={(n-k)!\over(n-2k)!}P_{n-2k}$$
for $0\leq 2k\leq n$. This proves  \link{$(11)$}{e11}.\par
{$2^\rond.$} Just combine \link{$(11)$}{e11} and \link{$(9)$}{e9}.\cqfd
\bg 
\bookm{2}{Determinant identities}
\parindent=20pt
\proclaim Proposition 6.  The sequences $(P_n)_{n\in\nat}$  and $(Q_n)_{n\in\nat}$ satisfy the following properties.
\item{$1^\rond.$}  For all $n\geq0$, we have $Q_{n+1}P_n-P_{n+1}Q_n=(-1)^nn!$.\dest{e13}\hfill$(13)$
\item{$2^\rond.$}  For all $n\geq0$, we have $Q_{n+2}P_n-P_{n+2}Q_n=(-1)^nn!X$.\dest{e14}\hfill$(14)$
\par
\parindent=0pt
{\bf Proof.}\par
{$1^\rond.$} First, let
 $$\delta_n^{(1)}=\det\left[\matrix{Q_{n+1}&P_{n+1}\cr Q_n&P_n}\right]=Q_{n+1}P_n-P_{n+1}Q_n$$ 
Using \link{$(4)$}{e4} and  \link{$(5)$}{e5} we can write
$$ \left[\matrix{Q_{n+1}&P_{n+1}\cr Q_n&P_n}\right]=
\left[\matrix{X&n\cr 1&0}\right]\left[\matrix{Q_{n}&P_{n}\cr Q_{n-1}&P_{n-1}}\right]$$
This proves that $\delta_n^{(1)}=-n\delta_{n-1}^{(1)}$, and, in its turn, proves \link{$(13)$}{e13} by induction since $\delta_0^{(1)}=1$.\par
{$2^\rond.$} Let $U_n$ denote $P_n$ or $Q_n$. Using \link{$(4)$}{e4} and \link{$(5)$}{e5} we have, for $n\geq2$,
$$\eqalign{U_{n+2}=& XU_{n+1}+(n+1)U_n= X(XU_{n}+nU_{n-1})+(n+1)U_n\cr
=& (X^2+n+1)U_{n}+n(XU_{n-1})= (X^2+n+1)U_{n}+n(U_{n}-(n-1)U_{n-2})\cr
=& (X^2+2n+1)U_{n}-n(n-1)U_{n-2}}$$
From this recurrence relation we conclude that
$$ \left[\matrix{Q_{n+2}&P_{n+2}\cr Q_n&P_n}\right]=
\left[\matrix{X^2+2n+1&-n(n-1)\cr 1&0}\right]\left[\matrix{Q_{n}&P_{n}\cr Q_{n-2}&P_{n-2}}\right]$$
So, if we define,
 $$\delta_n^{(2)}={1\over n!}\det\left[\matrix{Q_{n+2}&P_{n+2}\cr Q_n&P_n}\right]={1\over n!}(Q_{n+2}P_n-P_{n+2}Q_n)$$ 
we conclude from the matrix identity that $\forall\,n\geq2,~\delta_n^{(2)}=\delta_{n-2}^{(2)}$. 
This proves \link{$(14)$}{e14} by induction since $\delta_0^{(2)}=X$ and $\delta_1^{(2)}=-X$.\cqfd\par
\bg
\bookm{1}{First order inequalities}
\noindent{\bf 3. First order inequalities.}
\bg
Our first family of inequalities gives upper and lower bounds for $\varphi$ on $\real_+^*$ using rational functions.
\bg
\bookm{2}{Inequalities using rational functions}
\parindent=20pt
\proclaim Proposition 7. \dest{p7} The sequences $(P_n)_{n\in\nat}$  and $(Q_n)_{n\in\nat}$ satisfy the following properties.
\item{$1^\rond.$}  For all $n\geq0$, and all $x>0$, we have
 $${Q_{2n}(x)\over P_{2n}(x)}<\varphi(x)<{Q_{2n+1}(x)\over P_{2n+1}(x)}.\dest{e15}\eqno(15)$$
\item{$2^\rond.$}  For all $n\geq0$, and all $x>0$, we have
 $$\abs{\varphi(x)-{Q_n(x)\over P_n(x)}}<{n!\over P_n(x)P_{n+1}(x)}.\dest{e16}\eqno(16)$$
\par
\parindent=0pt
{\bf Proof.}\par
{$1^\rond.$} Using the fact that $\varphi(x)=\int_0^\infty e^{-x t-t^2/2}\d t$ for  $x\in \real$ we conclude that
$$\forall\,n\in\nat,~\forall\,x\in\real,\quad \varphi^{(n)}(x)=(-1)^n\int_0^\infty t^n e^{-x t-t^2/2}\d t$$
and this implies, in particular, that $\forall\,n\in\nat,~\forall\,x\in\real,\quad (-1)^n\varphi^{(n)}(x)>0$. \par
$\bullet$ The non-trivial part of the inequality  $\forall\,n\in\nat,~\forall\,x\in\real,~ \varphi^{(2n)}(x)>0$ is equivalent to
$$\forall\,n\in\nat,~\forall\,x\in\real,\quad {Q_{2n}(x)\over P_{2n}(x)}<\varphi(x)$$
$\bullet$ Similarly, the non-trivial part of the inequality  $\forall\,n\in\nat,~\forall\,x\in\real,~\varphi^{(2n+1)}(x)<0$ is equivalent to
$$\forall\,n\in\nat,~\forall\,x\in\real,\quad \varphi(x)<{Q_{2n+1}(x)\over P_{2n+1}(x)}.$$
{$2^\rond.$} For a non-negative integer $n$, and a positive real $x$ we have
$${Q_{2n}(x)\over P_{2n}(x)}<\varphi(x)<{Q_{2n+1}(x)\over P_{2n+1}(x)}\qquad\hbox{and}\qquad
{Q_{2n+2}(x)\over P_{2n+2}(x)}<\varphi(x)<{Q_{2n+1}(x)\over P_{2n+1}(x)}$$
and using \link{$(13)$}{e13}, we conclude that :
$$0<\varphi(x)-{Q_{2n}(x)\over P_{2n}(x)}<{Q_{2n+1}(x)\over P_{2n+1}(x)}-{Q_{2n}(x)\over P_{2n}(x)} 
={(2n)!\over P_{2n}(x) P_{2n+1}(x)}$$ 
and
$$0<{Q_{2n+1}(x)\over P_{2n+1}(x)}-\varphi(x)<{Q_{2n+1}(x)\over P_{2n+1}(x)}-{Q_{2n+2}(x)\over P_{2n+2}(x)} 
={(2n+1)!\over P_{2n+1}(x) P_{2n+2}(x)}.$$ 
This ends the proof.\cqfd
\bg
\bookm{2}{Examples}
{\bf Examples.} Here are the some of these inequalities:\par
$$\matrix{n=0,&\ds 0<\varphi (x)<\frac{1}{x}\cr
n=1,&\ds \frac{x}{x^2+1}<\varphi (x)<\frac{x^2+2}{x^3+3 x}\cr
n=2,&\ds \frac{x^3+5 x}{x^4+6 x^2+3}<\varphi (x)<\frac{x^4+9 x^2+8}{x^5+10
   x^3+15 x}}$$
\bg
In fact, these inequalities are excellent, and become more and more accurate for the values of $x$ that are large enough, say $x\geq1$. 
Moreover, it is straightforward, using \link{$(10)$}{e10}, to see that
$$\forall\,x\in\real_+^*,\quad\lim_{n\to\infty}{Q_n(x)\over P_n(x)}=\varphi(x)$$
but unfortunately, this convergence is extremely slow, for $x$ near zero.\par
\bg
\bookm{2}{Continued fractions expansion}
{\bf Remark.} In view of the recurrence relations  \link{$(4)$}{e4}  and  \link{$(5)$}{e5} and the above result, we expect that things could be arranged to obtain the continued fractions expansion of $\varphi(x)$ for positive values of $x$.\par
Indeed, let us define the sequence $(a_n)_{n\in\nat}$ as follows
$$a_{2n}=\prod_{k\leq2n\atop k\,\hbox{\fiverm odd}}k={(2n)!\over2^nn!}\quad\hbox{and}\quad
a_{2n+1}=\prod_{k\leq2n+1\atop k\,\hbox{\fiverm even}}k=2^nn!$$
then, define $(b_n)_{n\in\nat}$ by $b_n=a_n/a_{n+1}$, so that
$$b_{2n}={(2n)!\over2^{2n}(n!)^2}={1\over2^{2n}}\left({2n\atop n}\right)\quad\hbox{and}\quad
b_{2n+1}={1\over (2n+1)b_{2n}}.$$
Finally, define $\ds \tilde{P}_n={1\over a_n}P_n$ and  $\ds \tilde{Q}_n={1\over a_n}Q_n$. Having set  this notation, it is now a simple matter to verify that
$$\matrix{(\tilde{P}_0,\tilde{P}_1)=(1,X),& \quad& \forall\,n\geq1,~\tilde{P}_{n+1}=b_nX\tilde{P}_n+\tilde{P}_{n-1}\cr\cr
(\tilde{Q}_0,\tilde{Q}_1)=(0,1),& \quad& \forall \,n\geq1,~\tilde{Q}_{n+1}=b_nX\tilde{Q}_n+\tilde{Q}_{n-1}}$$
and this proves that
$${Q_n\over P_n}={\tilde{Q}_n\over \tilde{P}_n}
=[0,b_0X,b_1X,\ldots,b_{n-1}X]$$
and gives the following continued fractions expansion of $\varphi(x)$ for positive values of $x$ :
$$\varphi(x)=[0,x,x,{x\over2},{2x\over3},{3x\over8},{8x\over15},\cdots,b_nx,\cdots].$$\par
In fact, a more beautiful formula is the following :
$$\forall\,x>0,\quad\varphi(x)={1\over{x+{1\over{x+{2\over{x+{3\over{x+{4\over{x+{5\over{x+{\scriptscriptstyle \cdots}}}}}}}}}}}}}$$\par
\bg
\eject
\bookm{1}{Second order inequalities}
\noindent{\bf 4. Second order inequalities.}
\bg
Our second family of inequalities gives upper and lower bounds for $\varphi$ using rational functions and simple square roots.
\bg
\bookm{2}{Log-convexity}
\parindent=20pt
\proclaim Proposition 8. For every non-negative integer $n$, the function $x\mapsto\ln\abs{\varphi^{(n)}(x)}$ is strictly convex on $\real$.
\par
\parindent=0pt
{\bf Proof.}\par
We have seen that
$$\forall\,n\in\nat,~\forall\,x\in\real,\quad \varphi^{(n)}(x)=(-1)^n\int_0^\infty t^n e^{-x t-t^2/2}\d t$$
so that $\abs{\varphi^{(n)}}=(-1)^n\varphi^{(n)}$. Now, for every real $\lambda$, we have
$$(-1)^n(\lambda^2\varphi^{(n)}(x)+2\lambda\varphi^{(n+1)}(x)+\varphi^{(n+2)}(x))=
\int_0^\infty (\lambda-t)^2t^n e^{-x t-t^2/2}\d t>0$$
this proves that the discriminant of the second degree equation $\lambda^2\varphi^{(n)}(x)+2\lambda\varphi^{(n+1)}(x)+\varphi^{(n+2)}(x)=0$, 
 with respect to the unknown $\lambda$, is negative, that is,
$$\forall\,n\in\nat,~\forall\,x\in\real,\quad\varphi^{(n)}(x)\varphi^{(n+2)}(x)-\left(\varphi^{(n+1)}(x)\right)^2
>0\dest{e17}\eqno(17)$$
and this inequality is equivalent to the fact that  $x\mapsto\ln\abs{\varphi^{(n)}(x)}$ is strictly convex on $\real$.\cqfd
\bg
As a corollary, we can prove inequality \link{$(18)$}{e18}, which is a well known inequalty of {\sc Komatsu} [\link{1}{r1}], 
and  inequality \link{$(19)$}{e19}, which was proved differently in [\link{2}{r2}].\bg
\bg
\bookm{2}{Known inequalities}
\parindent=20pt
\proclaim Corollary 9.  
\item{$1^\rond.$}  For all $x\in\real$, we have
 $\ds{2\over x+\sqrt{x^2+4}}<\varphi(x)$.\dest{e18}\hfill$(18)$
\item{$2^\rond.$}  For all $x>-1$, we have
 $\ds\varphi(x)<{4\over 3x+\sqrt{x^2+8}}$.\dest{e19}\hfill$(19)$
\par
\parindent=0pt
{\bf Proof.}\par
{$1^\rond.$} Clearly, $\varphi^\prime(x)=x\varphi(x)-1$ and
 $\varphi^{\prime\prime}(x)=(x^2+1)\varphi(x)-x$, so that, for $n=0$, inequality \link{$(17)$}{e17} is equivalent to
$$\forall\,x\in\real,\quad\varphi^2(x)+x\varphi(x)-1>0$$
Now, since $\varphi(x)$ is positive it must by larger than the positive root of the equation $T^2+xT-1=0$,
 with respect to the unknown $T$, that is
$$\forall\,x\in\real,\quad\varphi(x)>\sqrt{{x^2\over4}+1}-{x\over2}={2\over x+\sqrt{x^2+4}}$$
{$2^\rond.$} Again, noting that $\varphi^{\prime\prime\prime}(x)=(x^3+3x)\varphi(x)-x^2-2$, we deduce from inequality \link{$(17)$}{e17} 
for $n=1$, that
$$\forall\,x\in\real,\quad(x^2-1)\varphi^2(x)-3x\varphi(x)+2>0\dest{clb}\eqno(\clubsuit)$$
\parindent=20pt
\item{$\bullet$} For $x<-1$, inequality \link{$(\clubsuit)$}{clb} is trivial, since both roots of the equation $(x^2-1)T^2-3xT+2=0$, with respect to the unknown $T$, are negative in this case.
\item{$\bullet$} For $x>-1$, inequality \link{$(\clubsuit)$}{clb} is equivalent to
$$\left({4\over3x+\sqrt{x^2+8}}-\varphi(x)\right)\left({3x+\sqrt{x^2+8}\over2}-(x^2-1)\varphi(x)\right)>0$$
so that, the function $x\mapsto {4\over3x+\sqrt{x^2+8}}-\varphi(x)$, which is 
continuous, and does not vanish on the interval $]-1,+\infty[$, must have a
constant sign on this interval. On the other hand, it is positive for $x=1$, as it can be seen 
by substituting $x=1$ in \link{$(\clubsuit)$}{clb}. Consequently, we have proved that
$$\forall\,x>-1,\quad\varphi(x)<{4\over3x+\sqrt{x^2+8}}$$
\par
\parindent=0pt
Which is the desired result.\cqfd
\bg
We want to generalize the preceding inequalities, and to this end, we need the following technical lemma, whose proof is postponed to the end of this note.
\bookm{2}{A new polynomial sequence An}
\parindent=20pt
\proclaim Lemma 10. \dest{p10} Define $A_n=P_nP_{n+2}-(P_{n+1})^2$ for $n\in\nat$. Then, 
$$\eqalign{A_{2n}(X)&=(2n+1)\prod_{k=1}^n(2k-1)^2+\sum_{m=2}^{2n}\left({(2n)!\over m!}\sum_{k=0}^{n-\ceiling{m/2}}\left({2n-2k-2\atop m-2}\right){(2k+1)!\over 2^{2k}(k!)^2}\right)X^{2m}\cr
A_{2n+1}(X)&=\left(\prod_{k=1}^{n+1}(2k-1)^2\right)(X^2-1)+\sum_{m=2}^{2n+1}\left({(2n+1)!\over m!}\sum_{k=0}^{n-\floor{m/2}}\left({2n-2k-1\atop m-2}\right){(2k+1)!\over 2^{2k}(k!)^2}\right)X^{2m}}$$
In particular,
\item{$\bullet$} $A_{2n}$ is an even polynomial of degree $4n$, whose coefficients are non-negative integers, and particularly, it does not have real roots.
\item{$\bullet$} $A_{2n+1}-A_{2n+1}(0)$ is an even polynomial of degree $4n+2$, whose coefficients are non-negative integers, and particularly, there is a unique $\beta_n$ from the interval $]0,1]$, such that
$$\eqalign{\abs{x}<\beta_n&\Rightarrow A_{2n+1}(x)<0\cr
\abs{x}>\beta_n&\Rightarrow A_{2n+1}(x)>0}$$\par
\bg
\parindent=0pt
We come now to the main theorem, which contains the second family of inequalities.\bg 
\bookm{2}{The second family of inequalities}
\parindent=20pt
\proclaim Theorem 11. \dest{p11} Consider the two sequences of polynpmials $(A_n)_{n\in\nat}$ and $(B_n)_{n\in\nat}$ defined by
$$A_n=P_nP_{n+2}-(P_{n+1})^2\qquad\hbox{and}\qquad B_n=P_nQ_{n+2}+P_{n+2}Q_n-2P_{n+1}Q_{n+1}.$$
Then, for every non-negative integer $m$, we have 
$$\forall\, x \in \real,\quad \varphi(x)>{B_{2m}(x)+(2m)!\,\sqrt{x^2+8m+4}\over 2A_{2m}(x)}\dest{e20}\eqno{({\cal I}_{2m})}$$
$$\forall\, x \in ]-\beta_m,+\infty[,\quad \varphi(x)<{B_{2m+1}(x)-(2m+1)!\,\sqrt{x^2+8m+8}\over 2A_{2m+1}(x)}\dest{e21}
\eqno{({\cal I}_{2m+1})}$$
Moreover, these inequalities are more precise than those in \link{proposition 7}{p7}.\par
\parindent=0pt
{\bf Proof.}\par
The inequality of log-convexity \link{$(17)$}{e17} is equivalent, for $n\in\nat$ and $x\in\real$, to the following inequality, 
$$(P_n(x)\varphi(x)-Q_n(x))(P_{n+2}(x)\varphi(x)-Q_{n+2}(x))-(P_{n+1}(x)\varphi(x)-Q_{n+1}(x))^2>0$$
So, if we define also $C_n=Q_nQ_{n+2}-(Q_{n+1})^2$, we see that,  for $n\in\nat$ and $x\in\real$, we have
$$A_n(x)(\varphi(x))^2-B_n(x)\varphi(x)+C_n(x)>0\dest{hrt}\eqno(\heartsuit)$$
-- The surprise here, is the very simple expression of $\Delta_n(X)$ the discriminant of the second degree equation
$$A_n(X)T^2-B_n(X)T+C_n(X)=0$$
with respect to the unknown $T$. In fact, by applying the identity
$$(ac^\prime+ca^\prime-2bb^\prime)^2-4(ac-b^2)(a^\prime c^\prime-b^{\prime2})
=(ac^\prime-ca^\prime)^2-4(ab^\prime-ba^\prime)(bc^\prime-cb^\prime)$$
to $(a,b,c)=(P_n,P_{n+1},P_{n+2})$ and $(a^\prime,b^\prime,c^\prime)=(Q_n,Q_{n+1},Q_{n+2})$, we conclude that
$$\Delta_n=(P_{n}Q_{n+2}-P_{n+2}Q_{n})^2-(P_{n}Q_{n+1}-P_{n+1}Q_{n})(P_{n+1}Q_{n+2}-P_{n+2}Q_{n+1})$$
and using \link{$(13)$}{e13} and \link{$(14)$}{e14} we come to the conclusion that
$$\Delta_n(X)=(n!)^2X^2-4(-1)^nn!\times(-1)^{n+1}(n+1)!=(n!)^2(X^2+4n+4).$$
-- That is, the equation $A_n(x)T^2-B_n(x)T+C_n(x)=0$, with respect to the unknown $T$,  has two different real roots $Z_n^+(x)$ and
 $Z_n^-(x)$ given by
$$Z_n^\pm(x)={B_n(x)\pm n!\,\sqrt{x^2+4n+4}\over 2A_n(x)}$$
(except for $n=2m+1$ and $x=\beta_m$, see \link{lemma 10}{p10}).\par
-- In order to determine the relative position of the approximation $Q_n(x)/P_n(x)$ of $\varphi(x)$ with respect to the roots  $Z_n^+(x)$ and
 $Z_n^-(x)$, we evaluate of the trinomial (in $T$) $A_n(x)T^2-B_n(x)T+C_n(x)$ for $T=Q_n(x)/P_n(x)$. Using some algebra and making use of \link{$(13)$}{e13}, we find 
$$A_n\left({Q_n\over P_n}\right)^2-B_n\left({Q_n\over P_n}\right)+C_n
=-{(P_{n}Q_{n+1}-P_{n+1}Q_{n})^2\over P_n^2}=-\left({n!\over P_n}\right)^2$$
This means that $Q_n(x)/P_n(x)$ falls between the two roots  $Z_n^+(x)$ and
 $Z_n^-(x)$ in the case $A_n(x)>0$ and falls outside the interval defined by $Z_n^+(x)$ and
 $Z_n^-(x)$ if $A_n(x)<0$.\par
-- {\bf The case} $n=2m$.\par
By \link{lemma 10}{p10} we have $\forall\, x \in \real,~A_{2m}(x)>0$, so that $Z_{2m}^-(x)<Z_{2m}^+(x)$, 
and \link{$(\heartsuit)$}{hrt} proves that $\varphi(x)$ falls outside the interval  $[Z_{2m}^-(x),Z_{2m}^+(x)]$, but using 
\link{$(15)$}{e15} and the above point,  we find
$$\forall\, x>0,\quad Z_{2m}^-(x)<{Q_{2m}(x)\over P_{2m}(x)}<\varphi(x)$$
so that $\forall\, x>0,\quad Z_{2m}^+(x)<\varphi(x)$.\par
The continuous function defined on $\real$ by
 $x\mapsto \varphi(x)-Z_{2m}^+(x)$ does not vanish on  $\real$ according to \link{$(\heartsuit)$}{hrt}, 
consequently, it has a constant sign, but it is positive for $x>0$, so it must be positive on $\real$, and this proves \link{$({\cal I}_{2m})$}{e20}. On the 
other hand, since $Q_{2m}(x)/P_{2m}(x)$ belongs to  $]Z_{2m}^-(x),Z_{2m}^+(x)[$ for $x>0$ we conclude that
$$\forall\, x \in \real_+^*,\quad {Q_{2m}(x)\over P_{2m}(x)}<Z_{2m}^+(x)<\varphi(x)$$
so that, inequality \link{$({\cal I}_{2m})$}{e20} is more precise than \link{$(15)$}{e15}.\par
-- {\bf The case} $n=2m+1$. Here we consider two possibilities :
\parindent=20pt
\item{$\bullet$} By \link{lemma 10}{p10} we have $\forall\, x >\beta_m,~A_{2m+1}(x)>0$, so that $Z_{2m+1}^-(x)<Z_{2m+1}^+(x)$, 
and \link{$(\heartsuit)$}{hrt} proves that $\varphi(x)$ falls outside the interval  $[Z_{2m+1}^-(x),Z_{2m+1}^+(x)]$, but using 
\link{$(15)$}{e15} and the above point,  we find
$$\forall\, x>\beta_m,\quad \varphi(x)<{Q_{2m+1}(x)\over P_{2m+1}(x)}<Z_{2m+1}^+(x)$$
so that $$\forall\, x>\beta_m,\quad \varphi(x)<Z_{2m+1}^-(x)<{Q_{2m+1}(x)\over P_{2m+1}(x)}$$
this proves \link{$({\cal I}_{2m+1})$}{e21}, for $x>\beta_m$, and proves that it is more precise than \link{$(15)$}{e15} in this case.
\item{$\bullet$} On the other hand, again by \link{lemma 10}{p10} we have $\forall\, x \in\,]-\beta_m,\beta_m[,~A_{2m+1}(x)<0$, so that $Z_{2m+1}^+(x)<Z_{2m+1}^-(x)$, 
and it follows from \link{$(\heartsuit)$}{hrt} that   $Z_{2m+1}^+(x)<\varphi(x)<Z_{2m+1}^-(x)$. This proves inequality \link{$({\cal I}_{2m+1})$}{e21}, for $-\beta_m<x<\beta_m$. We also have
$$\forall\, x\in ]0,\beta_m,[\quad \varphi(x)<Z_{2m+1}^-(x)<{Q_{2m+1}(x)\over P_{2m+1}(x)}$$
because $Q_{2m+1}(x)/ P_{2m+1}(x)$ falls outside the interval $[Z_{2m+1}^+(x),Z_{2m+1}^-(x)]$
for $0<x<\beta_m$.
\item{$\bullet$} Finally, the case $x=\beta_m$ can be treated by a simple continuity argument.
\par
\parindent=0pt
This finishes the proof of theorem 11.\cqfd
\bg
\bookm{2}{Examples}
{\bf Examples.} The two inequalities that come just after \link{$(18)$}{e18} and  \link{$(19)$}{e19} are 
$$\matrix{{\cal I}_2~:&\ds\forall\, x \in \real,\quad{x^3-2x+\sqrt{x^2+12}\over x^4+3}<\varphi(x)\cr
{\cal I}_3~:&\ds\forall\, x \in\, ]-\beta_1,+\infty[,\quad \varphi(x)<{x^4+x^2+16\over x^5+2x^3+12x+3\sqrt{x^2+16}}}$$
with $\beta_1\approx 0.871\,338$.\bg
{\bf Remark.} In the proof of \link{theorem 11}{p11}, we did not need the explicit form of $A_n(X)$ given in \link{lemma 10}{p10}, but expliciting $A_n$  was the only way, we could find, to study the variations of these polynomials. So, we will end this note by presenting a proof of this lemma.\bg
\bookm{2}{Proof of lemma 10}
{\bf Proof of lemma 10.}\par
In order to determine $A_n$ we will use generating functions. Let us start by finding an integral formula
for $A_n$. Using \link{$(8)$}{e8} we can write
$$\eqalign{P_{n}(x)P_{n+2}(x)&={1\over2\pi}\int_{\real^2} t^nu^{n+2}\exp\left(-{(t-x)^2+(u-x)^2\over2}\right)\d t\d u\cr
&={1\over2\pi}\int_{\real^2}t^{n+2}u^{n}\exp\left(-{(t-x)^2+(u-x)^2\over2}\right)\d t\d u\cr
P_{n+1}(x)P_{n+1}(x)&={1\over2\pi}\int_{\real^2} t^{n+1}u^{n+1}\exp\left(-{(t-x)^2+(u-x)^2\over2}\right)\d t\d u}$$
and from $A_n=P_{n}P_{n+2}-P_{n+1}^2$ we conclude that
$$A_n(x)={1\over4\pi}\int_{\real^2}( tu)^n(u-t)^2\exp\left(-{(t-x)^2+(u-x)^2\over2}\right)\d t\d u\dest{spd}\eqno(\spadesuit)$$
-- Consider $\rho\in [0,1[$ and $n\in\nat$, then from \link{$(\spadesuit)$}{spd} we can write
$$\eqalign{{\abs{A_n(x)}\rho^n\over n!}&\leq{1\over4\pi}\int_{\real^2}{\abs{ tu}^n\rho^n\over n!}(u-t)^2\exp\left(-{(t-x)^2+(u-x)^2\over2}\right)\d t\d u\cr
&\leq{1\over4\pi}\int_{\real^2}(u-t)^2\exp(\rho\abs{ tu})\exp\left(-{(t-x)^2+(u-x)^2\over2}\right)\d t\d u\cr
&\leq{1\over4\pi}\int_{\real^2}(u-t)^2\exp\left(\rho{t^2+u^2\over2}\right)\exp\left(-{(t-x)^2+(u-x)^2\over2}\right)\d t\d u\cr
&\leq{e^{-x^2}\over4\pi}\int_{\real^2}(u-t)^2\exp\left(x(t+u)-(1-\rho){t^2+u^2\over2}\right)\d t\d u.
}$$
and since $0\leq\rho<1$ we conclude that $\ds\sup_{n\geq0}\left({\abs{A_n(x)}\rho^n\over n!}\right)<+\infty$. This proves that the
radius of convergence of the series $\ds y\mapsto \sum_{n=0}^\infty{A_n(x)\over n!}y^n$ is at least 1.\par
-- Let us then consider, for $x\in \real$, the analytic function $G_x$ defined by
$$G_x:]-1,1[\rightarrow\real, G_x(y)=\sum_{n=0}^\infty{A_n(x)\over n!}y^n$$
using \link{$(\spadesuit)$}{spd}, we have
$$\eqalign{G_x(y)&={1\over4\pi}\int_{\real^2}\left(\sum_{n=0}^\infty{(ut)^ny^n\over n!}\right)(t-u)^2\exp\left(-{(t-x)^2+(u-x)^2\over2}\right)\d t\d u\cr
&={1\over4\pi}\int_{\real^2}(t-u)^2\exp\left(uty-{(t-x)^2+(u-x)^2\over2}\right)\d t\d u}$$
Now, using the change of variables $u\leftarrow v-w$ and $t\leftarrow v+w$ we obtain
$$G_x(y)={2\over \pi}\underbrace{\left(\int_\real e^{v^2y-(x-v)^2}\d v\right)}_{\cal J}
\underbrace{\left(\int_\real w^2 e^{-(1+y)w^2}\d w\right)}_{\cal K}$$
\parindent=20pt
\item{$\bullet$} To calculate ${\cal K}$ we use the change of variable  $\ds w\leftarrow {s\over\sqrt2(1+y)}$ :
$${\cal K}={1\over2\sqrt2(1+y)\sqrt{1+y}}\int_\real s^2 e^{-s^2/2}\d s={\sqrt\pi\over2(1+y)\sqrt{1+y}}$$
\item{$\bullet$} As for ${\cal J}$, we note that $(x-v)^2-v^2y=(1-y)\left(v-{x\over1-y}\right)^2-{yx^2\over1-y}$, which allows us to write  :
$$\eqalign{{\cal J}&=\exp\left({yx^2\over1-y}\right)\int_\real\exp\left(-(1-y)\left(v-{x\over1-y}\right)^2\right)\d v\cr
&=\exp\left({yx^2\over1-y}\right)\int_\real\exp\left(-(1-y)v^2\right)\d v\cr
&=\exp\left({yx^2\over1-y}\right){1\over\sqrt2\sqrt{1-y}}\int_\real e^{-s^2/2}\d s\cr
&={\sqrt\pi\over\sqrt{1-y}}\exp\left({yx^2\over1-y}\right)}$$\par
\parindent=0pt
So, we have proved that $G_x$ has the following expression
$$G_x(y)={1\over(1+y)\sqrt{1-y^2}}\exp\left({yx^2\over1-y}\right)$$
and by expanding the exponential function we obtain
$$G_x(y)={1-y\over(1-y^2)^{3/2}}+{y\over(1-y^2)^{3/2}}x^2+\sum_{m=2}^\infty{y^m\over(1-y)^{m-1}(1-y^2)^{3/2}}\cdot {x^{2m}\over m!}.$$
Finally, using the expansions 
$$ {1\over (1-y^2)^{3/2}}=\sum_{k=0}^\infty{(2k+1)!\over2^{2k}(k!)^2}y^{2k}$$
and
$$ {y^m\over (1-y)^{m-1}}=\sum_{k=m}^\infty\left({k-2\atop m-2}\right)y^{k}$$
we deduce that
$$\eqalign{{1-y\over (1-y^2)^{3/2}}&=\sum_{n=0}^\infty{(-1)^n(n+1)\over 2^n}\left({n \atop \floor{n/2}}\right)y^{n}\cr
{y\over (1-y^2)^{3/2}}&=\sum_{n=0}^\infty{1-(-1)^n\over 2^n}n\left({n-1 \atop \floor{n/2}}\right)y^{n}\cr
{y^m\over (1-y)^{m-1}(1-y^2)^{3/2}}&=\sum_{n=m}^\infty\left(\sum_{k=0}^{\floor{(n-m)/2}}
{(2k+1)!\over 2^{2k}(k!)^2}\left({n-2k-2 \atop m-2}\right)\right)y^{n}}$$
So, if we define
$$a_{m,n}=\left\{\matrix{\ds{(-1)^n(n+1)\over 2^n}\left({n \atop \floor{n/2}}\right)&\hbox{if}&m=0\cr
\ds{1-(-1)^n\over 2^n}n\left({n-1 \atop \floor{n/2}}\right)&\hbox{if}&m=1\cr\ds
\sum_{k=0}^{\floor{(n-m)/2}}
{(2k+1)!\over 2^{2k}(k!)^2}\left({n-2k-2 \atop m-2}\right)&\hbox{if}&m\geq2}\right.$$
we get,
$$G_x(y)=\sum_{m=0}^\infty\left(\sum_{n=m}^\infty a_{n,m}y^n\right){x^{2m}\over m!}
=\sum_{n=0}^\infty\left(\sum_{m=0}^n{a_{n,m}\over m!}x^{2m}\right)y^n$$
and from this we obtain $\ds A_n(x)=\sum_{m=0}^n{a_{n,m}\over m!}x^{2m}$, which yields the desired expression.\cqfd
\par

\bg
\bookm{1}{References}
\noindent{\bf References.}\par
\parindent=20pt
\ref 1\dest{r1}& K.~Ito, H.~P.~McKean&Diffusion processes and their sample paths,& Springer-Verlag, & 1965. pp 17.
\ref 2\dest{r2}& S.~J.~Szarek, E.~Werner& &Confidence regions for means of multivariate normal distributions and
a non-symmetric correlation inequality for gaussian measure,& arXiv:math.FA/9701205.
\par
\parindent=0pt
\end
\item{}Department of Mathematics, {\sl  Higher Institute for Applied Sciences and Technology,} P.O. Box 31983, Damascus, Syria.
\item{} {\it E-mail} : omran\_kouba@hiast.edu.sy\par